# CHARACTERIZATION OF THE CUBIC EXPONENTIAL FAMILIES BY ORTHOGONALITY OF POLYNOMIALS

By Abdelhamid Hassairi and Mohammed Zarai

*Sfax University*

This paper introduces a notion of 2-orthogonality for a sequence of polynomials to give extended versions of the Meixner and Feinsilver characterization results based on orthogonal polynomials. These new versions subsume the Letac–Mora characterization of the real natural exponential families having cubic variance function.

**1. Introduction.** Let $F = \{P(m, F), m \in M_F\}$ be a natural exponential family (NEF) on the real line parameterized by its domain of the means $M_F$. If $V_F(m)$ denotes the variance of the probability distribution $P(m, F)$, then the mapping $m \mapsto V_F(m)$ is called the variance function of the family $F$. The importance of the variance function stems from the fact that it characterizes the family $F$ within the class of all natural exponential families. Furthermore, for many common NEFs the variance function takes a very simple form. Morris (1982) describes the class of real NEFs such that the variance function is a polynomial function of degree at most 2 in the mean. Up to affine transformations and powers of convolution, this class includes the normal, Poisson, binomial, negative binomial, gamma and a sixth family called hyperbolic cosine and nothing else. The Morris class of quadratic NEFs has received a deal of attention in the statistical literature and many interesting characteristic properties have been established. We will be concerned here only with the properties based on the notion of orthogonal polynomials. A remarkable characteristic result is due to Meixner (1934) [see also Letac (1992)]. It characterizes the distributions $\mu$ for which there exists a family of $\mu$-orthogonal polynomials with an exponential generating function. These distributions generate exactly the Morris class of NEFs. A second characterization is due to Feinsilver (1986), who shows that a certain class of polynomials naturally associated to a NEF is $\mu$-orthogonal if and only if the family is in the Morris class. In the present paper, we will be concerned with the class of real natural exponential families having cubic variance function (i.e., a variance function which is a polynomial of degree less than or equal to 3). In fact, Letac and Mora (1990) have extended the work of Morris by classifying all real cubic natural exponential families. They have added to the Morris class six other types of NEFs which may be obtained from the Morris class by the action of the linear group $GL(\mathbb{R}^2)$ [see Hassairi (1992)]. The best known among such families is the Inverse-Gaussian family. Our aim is then to extend to the Letac–Mora class of cubic NEFs the different characterizations established for the Morris class and based on orthogonal polynomials. We show in fact that the cubicity of the variance function is characterized by a property of orthogonality which will be called the property of 2-orthogonality. The result is interesting in its own right and seems potentially very useful for asymptotic expansions and approximation. In Section 2, after a review of exponential family theory, we specify some facts about the Feinsilver sequence of polynomials associated to a NEF; in particular, we show that the generating function of this sequence converges in a neighborhood of 0. In Section 3, we state and prove our main result concerning the characterization of the Feinsilver sequence of polynomials corresponding to a distribution generating a cubic natural exponential family by a property similar to orthogonality. We also show that this sequence is characterized by a four-term recurrence relation, while, as it is well known, a sequence of orthogonal polynomials satisfies a three-term recurrence relation. In Section 4, we first determine the families of 2-orthogonal polynomials with exponential generating function. [The families with a formal exponential generating function are sometimes called in the literature Sheffer polynomials of type 0; see Sheffer (1939) and Rainville (1960).] This leads to another characterization of the Letac–Mora class of cubic NEFs which may be considered as the extension to this class of the Meixner characterization.

**2. Exponential families and associated polynomials.** We need first to review some facts concerning natural exponential families and to introduce some notations. If $\mu$ is a positive Radon measure on the real line $\mathbb{R}$, we denote by

$$L_\mu(\theta) = \int_\mathbb{R} \exp(\theta x)\mu(dx) \leq +\infty$$

its Laplace transform, and we denote by $\Theta(\mu)$ the interior of the convex set $D(\mu) = \{\theta \in \mathbb{R}; L_\mu(\theta) < \infty\}$. $\mathcal{M}(\mathbb{R})$ will denote the set of measures $\mu$ such that $\Theta(\mu)$ is not empty and $\mu$ is not concentrated on one point. If $\mu$ is in $\mathcal{M}(\mathbb{R})$, we also denote

$$k_\mu(\theta) = \log L_\mu(\theta), \qquad \theta \in \Theta(\mu),$$

the cumulate function of $\mu$.

To each $\mu$ in $\mathcal{M}(\mathbb{R})$ and $\theta$ in $\Theta(\mu)$, we associate the following probability distribution on $\mathbb{R}$:

$$(2.1) \qquad P(\theta, \mu)(dx) = \exp(\theta x - k_\mu(\theta))\mu(dx).$$

The set

$$F = F(\mu) = \{P(\theta, \mu); \theta \in \Theta(\mu)\}$$

is called the NEF generated by $\mu$. If $\mu$ and $\mu'$ are in $\mathcal{M}(\mathbb{R})$, then $F(\mu) = F(\mu')$ if and only if there exists $(a, b)$ in $\mathbb{R}^2$ such that $\mu'(dx) = \exp(ax + b)\mu(dx)$. Therefore, if $\mu$ is in $\mathcal{M}(\mathbb{R})$ and $F = F(\mu)$,

$$\mathcal{B}_F = \{\mu' \in \mathcal{M}(\mathbb{R}); F(\mu') = F\}$$

is the set of basis of $F$.

We have, for all $\theta$ in $\Theta(\mu)$,

$$(2.2) \qquad \frac{d^n}{d\theta^n} L_\mu(\theta) = \int_\mathbb{R} x^n \exp(\theta x)\mu(dx) < +\infty,$$

and consequently, all the moments of $P(\theta, \mu)$ are finite. In fact,

$$(2.3) \qquad \frac{1}{L_\mu(\theta)} \frac{d^n}{d\theta^n} L_\mu(\theta) = \int_\mathbb{R} x^n P(\theta, \mu)(dx) < +\infty.$$

The function $k_\mu$ is strictly convex and real analytic. Its first derivative $k'_\mu$ defines a diffeomorphism between $\Theta(\mu)$ and its image $M_F$. Since $k'_\mu(\theta) = \int_\mathbb{R} xP(\theta, \mu)(dx)$, $M_F$ is called the domain of the means of $F$. The inverse function of $k'_\mu$ is denoted by $\psi_\mu$ and, setting $P(m, F) = P(\psi(m), \mu)$, the probability of $F$ with mean $m$, we have

$$F = \{P(m, F); m \in M_F\},$$

which is the parameterization of $F$ by the mean.

The density of $P(m, F)$ with respect to $\mu$ is

$$(2.4) \qquad f_\mu(x, m) = \exp\{\psi_\mu(m)x - k_\mu(\psi_\mu(m))\}.$$

For $m$ in $M_F$, we denote

$$(2.5) \qquad V_F(m) = \int_\mathbb{R} (x - m)^2 P(m, F)(dx).$$

Then

$$(2.6) \qquad V_F(m) = k''_\mu(\psi_\mu(m)) = (\psi'_\mu(m))^{-1},$$

and, the map $m \mapsto V_F(m)$ is called the variance function of $F$. It entirely characterizes the NEF; that is, if $F$ and $F'$ are two NEF such that $V_F(m) = V_{F'}(m)$ on a nonempty open set included in $M_F \cap M_{F'}$, then $F = F'$.

Consider now a real natural exponential family $F$ and take $\mu = P(m_o, F)$ with $m_o$ fixed in $M_F$. The density $f_\mu(\cdot, m)$ of $P(m, F)$ with respect to $\mu$ is still given by (2.4) with $f_\mu(\cdot, m_o) \equiv 1$. It is easily verified by induction on $n$ in $\mathbb{N}$ that there exists a polynomial $P_n$ in $x$ of degree $n$ such that

$$(2.7) \qquad \frac{\partial^n}{\partial m^n} f_\mu(x, m) = P_n(x, m) f_\mu(x, m)$$

and

$$(2.8) \qquad P_{n+1}(x, m) = \psi'_\mu(m)(x - m) P_n(x, m) + R_{n+1}(x, m),$$

where $R_{n+1}$ is a polynomial in $x$ of degree $< n+1$. In particular, we have that

$$(2.9) \qquad P_o(x, m) = 1 \quad \text{and} \quad P_1(x, m) = \psi'_\mu(m)(x - m).$$

We now make a useful observation through the following theorem whose proof will be given in Section 5. For the sake of simplification, we set

$$(2.10) \qquad P_n(x) = P_n(x, m_o).$$

THEOREM 2.1. *Let $F$ be a NEF on $\mathbb{R}$ and let $\mu$ be a fixed probability in $F$ with mean $m_o$. Let $P_n(x)$ be the polynomials defined by (2.7). Then*

$$\sum_{n \in \mathbb{N}} \frac{(m - m_o)^n}{n!} P_n(x)$$

*is an entire series in $L^2(\mu)$ of nonzero radius of convergence.*

It should be remarked that there exists $r > 0$ such that, for all $m$ in $]m_o - r, m_o + r[$ and for all $x$ in $\mathbb{R}$,

$$(2.11) \qquad f_\mu(x, m) = \sum_{n \in \mathbb{N}} \frac{(m - m_o)^n}{n!} P_n(x).$$

To conclude this section, we mention that in many interesting situations, $P_n$ can be calculated by mean of the Faà di Bruno formula

$$(f_o g)^{(n)}(m) = \sum \frac{n!}{k_1! \cdots k_n!} f^{(k)}(g(m)) \left(\frac{g^{(1)}(m)}{1!}\right)^{k_1} \cdots \left(\frac{g^{(n)}(m)}{n!}\right)^{k_n},$$

where $k = k_1 + \cdots + k_n$ and the sum is taken for all integers $k_j \geq 0$ such that $k_1 + 2k_2 + \cdots + nk_n = n$. That is, if we denote

$$g^{(k)}(x) = \frac{\partial^k}{\partial m^k}[(\psi_\mu(m)x - k_\mu(\psi_\mu(m)))]_{m=m_o},$$

then

$$P_n(x) = \sum_{k_1 + 2k_2 + \cdots + nk_n = n} \frac{n!}{k_1! \cdots k_n!} \left(\frac{g^{(1)}(x)}{1!}\right)^{k_1} \cdots \left(\frac{g^{(n)}(x)}{n!}\right)^{k_n}$$

with $k_j$ in $\mathbb{N}$, for $1 \leq j \leq n$.

The most famous example in this topic is the inverse Gaussian distribution with parameters $1/2$ and $p > 0$ defined by

$$\mu(dx) = \frac{p}{\sqrt{2\pi}} x^{-3/2} \exp\left(-\frac{p^2}{2x}\right) \mathbb{1}_{]0,+\infty[}(x)\, dx.$$

The NEF generated by $\mu$ belongs to the Letac–Mora class. In fact, a standard computation shows that $\Theta(\mu) = ]-\infty, 0[$ and $k_\mu(\theta) = -p\sqrt{-2\theta}$. If $F = F(\mu)$, then $M_F = ]0, +\infty[$, $\psi_\mu(m) = -p^2/2m^2$ and $V_F(m) = m^3/p^2$.

For all $m_o \in M_F$, we have

$$P_n(x) = \sum_{k_1+2k_2+\cdots+nk_n=n} \frac{n!}{k_1!\cdots k_n!} p^{2n} \left(\frac{2}{m_o^3}x - \frac{1}{m_o^2}\right)^{k_1} \cdots \left(\frac{n+1}{m_o^{n+2}}x - \frac{1}{m_o^{n+1}}\right)^{k_n}.$$

**3. Characterization of the cubic families in the Feinsilver way.** As pointed out in the Introduction, any real cubic natural exponential family can be obtained from a quadratic family via the action of the linear group $GL(\mathbb{R}^2)$. For instance, the quadratic natural exponential families on $\mathbb{R}$ have been characterized by Feinsilver (1986) as the ones for which the polynomials $P_n(x)$ are $\mu$-orthogonal. In this section we show that the polynomials $P_n(x)$ associated to a cubic natural exponential family have also a characterizing property of orthogonality which will be called the property of 2-orthogonality.

DEFINITION 3.1. Let $\mu$ be a measure on $\mathbb{R}$ such that $\int |x|^n \mu(dx) < \infty$ for all $n \in \mathbb{N}$. A family $(Q_n)_{n \in \mathbb{N}}$ of polynomials on $\mathbb{R}$ is $\mu - 2$-orthogonal if, for all $n$ and $q$ in $\mathbb{N}^*$, $\int Q_n(x) Q_q(x) \mu(dx) = 0$ when $n \geq 2q$, and $\int Q_n(x) \mu(dx) = 0$.

Next we give our first main result.

THEOREM 3.1. *Let $F$ be a NEF on $\mathbb{R}$ and let $\mu$ be an element of $F$ with mean $m_o$. Consider the polynomials $(P_n)_{n \in \mathbb{N}}$ defined by $P_n(x) = \frac{\partial^n}{\partial m^n} f_\mu(x, m)|_{m=m_o}$. Then the three following statements are equivalent:*

(i) *The polynomials $(P_n)_{n \in \mathbb{N}}$ are $\mu - 2$-orthogonal.*
(ii) *$F$ is cubic.*
(iii) *There exist real numbers $(a_i)_{0 \leq i \leq 3}$ such that, for all $n \geq 2$,*

$$xP_n(x) = a_3 A_n^3 P_{n-2}(x) + n(a_2(n-1) + 1)P_{n-1}(x)$$
$$+ (na_1 + m_o)P_n(x) + a_o P_{n+1}(x),$$

*with $A_n^3 = n(n-1)(n-2)$. Furthermore, in this case we have*

$$V_F(m) = a_3(m - m_o)^3 + a_2(m - m_o)^2 + a_1(m - m_o) + a_o.$$

To help in the proof of this theorem, let us give in Table 1, for each of the six types of NEF on $\mathbb{R}$ with polynomial variance function of degree 3 [see Letac and Mora (1990)], the sequence of $P(m_o, F) - 2$-orthogonal polynomials $P_n(x)$ defined by its recurrence relation for $m_o = 1$.

PROOF OF THEOREM 3.1. (i) $\Rightarrow$ (ii). Equation (2.11) allows us to claim the existence of $r > 0$ such that, for all $m$ in the interval $]m_o - r, m_o + r[$ and for all $x$ in $\mathbb{R}$,

$$f_\mu(x, m) = \sum_{n \in \mathbb{N}} \frac{(m - m_o)^n}{n!} P_n(x).$$

If, for $(m, m') \in (]m_o - r, m_o + r[)^2$, we set

(3.1) $\quad g(m, m') = \exp\{k_\mu(\psi_\mu(m) + \psi_\mu(m')) - k_\mu(\psi_\mu(m)) - k_\mu(\psi_\mu(m'))\},$

then the $\mu - 2$-orthogonality of the polynomials $(P_n)$ and Theorem 2.1 imply that

$$g(m, m') = \int f_\mu(x, m) f_\mu(x, m') \mu(dx)$$

$$= \int \sum_{n,q \in \mathbb{N}} \frac{(m - m_o)^q (m' - m_o)^n}{n! q!} P_n(x) P_q(x) \mu(dx)$$

$$= 1 + \sum_{n,q \in \mathbb{N}^*} \frac{(m - m_o)^q (m' - m_o)^n}{n! q!} \int P_n(x) P_q(x) \mu(dx)$$

$$= 1 + \sum_{q \in \mathbb{N}^*, n \in [(q+1)/2, 2q-1]} \frac{(m - m_o)^q (m' - m_o)^n}{n! q!} \int P_n(x) P_q(x) \mu(dx).$$

Taking the derivative of (3.1) with respect to $m$, we get, for all $(m, m') \in (]m_o - r, m_o + r[)^2$,

(3.2) $\quad \psi'_\mu(m)(k'_\mu(\psi_\mu(m) + \psi_\mu(m')) - k'_\mu(\psi_\mu(m))) g(m, m')$
$$= \sum_{q \geq 1, n \in [(q+1)/2, 2q-1]} q\, a_{nq} (m - m_o)^{q-1} (m' - m_o)^n,$$

with $a_{nq} = \frac{1}{n! q!} \int P_n(x) P_q(x) \mu(dx)$.

Making $m = m_o$ in (3.2), then, since $\psi_\mu(m_o) = 0$, we get

$$\psi'_\mu(m_o)(m' - m_o) = a_{11}(m' - m_o).$$

This is true for all $m' \in ]m_o - r, m_o + r[$; then

$$a_{11} = \psi'_\mu(m_o).$$

Table 1

| Type | $m_o$ | Induction relations |
|---|---|---|
| Inverse Gaussian with parameter $p=1$ <br> $V_F(m) = (m-m_o)^3 + 3m_o(m-m_o)^2$ <br> $\quad + 3m_o^2(m-m_o) + m_o^3$ | 1 | $P_o(x) = 1$, <br> $P_1(x) = (x-1)$, <br> $P_2(x) = x^2 - 6x + 3$, <br> $P_{n+1}(x) = (x - 3n - 1)P_n(x)$ <br> $\quad - n(3n-2)P_{n-1}(x)$ <br> $\quad - A_n^3 P_{n-2}(x), \quad n \geq 2$ |
| Strict arcsine with parameter $p=1$ <br> $V_F(m) = (m-m_o)^3 + 3m_o(m-m_o)^2$ <br> $\quad + (3m_o^2 + 1)(m-m_o) + m_o^3 + m_o$ | 1 | $P_o(x) = 1$, <br> $P_1(x) = \frac{1}{2}(x-1)$, <br> $P_2(x) = \frac{1}{4}(x^2 - 5x + 3)$, <br> $P_{n+1}(x) = \frac{1}{2}[(x - 4n - 1)P_n(x)$ <br> $\quad - n(3n-2)P_{n-1}(x)$, <br> $\quad - A_n^3 P_{n-2}(x)]$, <br> $\quad n \geq 2$ |
| Takács with parameters $p=1$ and $a=1$ <br> $V_F(m) = 2(m-m_o)^3$ <br> $\quad + (6m_o + 3)(m-m_o)^2$ <br> $\quad + (6m_o^2 + 6m_o + 1)(m-m_o)$ <br> $\quad \times (2m_o^3 + 3m_o^2 + m_o)$ | 1 | $P_o(x) = 1$, <br> $P_1(x) = \frac{1}{6}(x-1)$, <br> $P_2(x) = \frac{1}{36}(x^2 - 15x + 8)$, <br> $P_{n+1}(x) = \frac{1}{6}[(x - 13n - 1)P_n(x)$ <br> $\quad - n(9n-8)P_{n-1}(x)$ <br> $\quad - 2A_n^3 P_{n-2}(x)]$, <br> $\quad n \geq 2$ |
| Large arcsine with parameters $p=1$ <br> $\quad$ and $a=1$ <br> $V_F(m) = 2(m-m_o)^3$ <br> $\quad + (6m_o + 2)(m-m_o)^2$ <br> $\quad + (6m_o^2 + 4m_o + 1)(m-m_o)$ <br> $\quad + 6m_o^3 + 2m_o^2 + m_o$ | 1 | $P_o(x) = 1$, <br> $P_1(x) = \frac{1}{9}(x-1)$, <br> $P_2(x) = \frac{1}{81}(x^2 - 13x + 3)$, <br> $P_{n+1}(x) = \frac{1}{9}[(x - 11n - 1)P_n(x)$ <br> $\quad - n(8n-7)P_{n-1}(x)$ <br> $\quad - 2A_n^3 P_{n-2}(x)]$, <br> $\quad n \geq 2$ |
| Ressel with parameter $p=1$ <br> $V_F(m) = 2(m-m_o)^3 + (6m_o + 3)(m-m_o)^2$ <br> $\quad + (6m_o^2 + 6m_o + 1)(m-m_o)$ <br> $\quad \times (2m_o^3 + 3m_o^2 + m_o)$ | 1 | $P_o(x) = 1$, <br> $P_1(x) = \frac{1}{2}(x-1)$, <br> $P_2(x) = \frac{1}{4}(x^2 - 7x + 4)$, <br> $P_{n+1}(x) = \frac{1}{2}[(x - 5n - 1)P_n(x)$ <br> $\quad - n(4n-3)P_{n-1}(x)$ <br> $\quad - A_n^3 P_{n-2}(x)]$, <br> $\quad n \geq 2$ |
| Abel with parameter $p=1$ <br> $V_F(m) = (m-m_o)^3 + (3m_o + 2)(m-m_o)^2$ <br> $\quad + (3m_o^2 + 4m_o + 1)(m-m_o)$ <br> $\quad \times (m_o^3 + 2m_o^2 + m_o)$ | 1 | $P_o(x) = 1$, <br> $P_1(x) = \frac{1}{4}(x-1)$, <br> $P_2(x) = \frac{1}{16}(x^2 - 15x + 8)$, <br> $P_{n+1}(x) = \frac{1}{4}[(x - 8n - 1)P_n(x)$ <br> $\quad - n(5n-4)P_{n-1}(x)$ <br> $\quad - A_n^3 P_{n-2}(x)]$, <br> $\quad n \geq 2$ |

Again we take the derivative of (3.2) with respect to $m$ and we let $m = m_o$. We get, for all $m'$ in $]m_o - r, m_o + r[$,

$$\psi_\mu''(m_o)(m' - m_o) + a_{11}^2(V_F(m') - V_F(m_o)) + a_{11}^2(m' - m_o)^2$$
$$= 2a_{22}(m' - m_o)^2 + 2a_{23}(m' - m_o)^3.$$

Therefore,

$$V_F(m') = 2a_{11}^{-2}a_{23}(m' - m_o)^3 + (2a_{22} - a_{11}^2)a_{11}^{-2}(m' - m_o)^2$$
$$- a_{11}^{-2}\psi_\mu''(m_o)(m' - m_o) + V_F(m_o).$$

This implies that $V_F$ is cubic on $]m_o - r, m_o + r[$ and, by extension, we obtain that $F$ is a cubic NEF.

(ii) $\Rightarrow$ (iii). From (ii), there exist real numbers $(a_i)_{0 \leq i \leq 3}$ such that

$$V_F = a_3(m - m_o)^3 + a_2(m - m_o)^2 + a_1(m - m_o) + a_o.$$

On the other hand, we know that there exists $r > 0$ such that, for all $m \in ]m_o - r, m_o + r[$ and for all $x \in \mathbb{R}$,

$$\sum_{n \in \mathbb{N}} \frac{(m - m_o)^n}{n!} P_n(x) = \exp\{\psi_\mu(m)x - k_\mu(\psi_\mu(m))\}.$$

Denoting $\theta = \psi_\mu(m)$, this may be written as

$$(3.3) \qquad \exp(\theta x) = \left(\sum_{n \in \mathbb{N}} (k_\mu'(\theta) - m_o)^n \frac{P_n(x)}{n!}\right) \exp\{k_\mu(\theta)\}.$$

Taking the derivative with respect to $\theta$ of (3.3) gives

$$x \exp(\theta x) = \sum_{n \in \mathbb{N}} \frac{P_n(x)}{n!}(n(m - m_o)^{n-1} k_\mu''(\theta) + (m - m_o)^n k'(\theta)) \exp\{k_\mu(\theta)\},$$

which is equivalent to

$$\sum_{n \in \mathbb{N}} \frac{(m - m_o)^n}{n!} x P_n(x)$$

$$= \sum_{n \in \mathbb{N}} \frac{P_n(x)}{n!}(n(m - m_o)^{n-1} k_\mu''(\theta) + (m - m_o)^n m)$$

$$= \sum_{n \in \mathbb{N}} \frac{P_n(x)}{n!}\left(n \sum_{k=0}^{3} a_k(m - m_o)^{n+k-1} + (m - m_o)^{n+1} + m_o(m - m_o)^n\right).$$

By identification, we get

$$xP_n(x) = a_3 n(n-1)(n-2) P_{n-2}(x)$$
$$+ n(a_2(n-1) + 1) P_{n-1}(x) + (na_1 + m_o) P_n(x) + a_o P_{n+1}(x),$$

and (iii) is proved.

(iii) $\Rightarrow$ (i). The result is easily obtained if we verify the three following facts:

(a) For all $n \in \mathbb{N}^*$, $\int P_n(x)\mu(dx) = 0$.
(b) There exist real numbers $\beta_{n,q}^s$ such that, for all $n, q \in \mathbb{N}^*$ verifying $n \geq 2q$,

$$x^q P_n(x) = \beta_{n,q}^0 P_{n-2q}(x) + \sum_{n-2q+1 \leq s \leq n+q} \beta_{n,q}^s P_s(x),$$

where $\beta_{n,q}^0 = 0$ if $n = 2q$ $\forall q \in \mathbb{N}^*$,

(c) There exist real numbers $(\alpha_q)_{0 \leq q \leq n}$ such that

$$P_n(x) = \alpha_n x^n + \sum_{0 \leq q \leq n-1} \alpha_q x^q.$$

*Proof of* (a). We first observe that

$$\int \frac{\partial}{\partial m} f_\mu(x, m) \mu(dx) = \psi'_\mu(m) \int (x - m) f_\mu(x, m) \mu(dx)$$

$$= \psi'_\mu(m) \int (x - m) P(m, F)(dx)$$

$$= 0.$$

Since, for all $n$, we have

$$\int \left| \frac{\partial^n}{\partial m^n} f_\mu(x, m) \right| \mu(dx) = \int |P_n(x - m)| f_\mu(x, m) \mu(dx)$$

$$= \int |P_n(x - m)| P(m, F)(dx) < +\infty,$$

[see (2.3) and (2.7)], then

$$\int \frac{\partial}{\partial m} \left\{ \frac{\partial^n}{\partial m^n} f_\mu(x, m) \right\} \mu(dx) = \frac{\partial}{\partial m} \int \frac{\partial^n}{\partial m^n} f_\mu(x, m) \mu(dx).$$

Hence we obtain that, for all $n \in \mathbb{N}^*$, $\int \frac{\partial^n}{\partial m^n} f_\mu(x, m) \mu(dx) = 0$. This, for $m = m_o$, gives $\int P_n(x) \mu(dx) = 0$.

*Proof of* (b). We can write (iii) as

(3.4) $$x P_n(x) = \beta_{n,1}^0 P_{n-2}(x) + \sum_{n-1 \leq s \leq n+1} \beta_{n,1}^s P_s(x),$$

where $\beta_{n,1}^0 = a_3 n(n-1)(n-2)$.

For a fixed $n$ in $\mathbb{N}^*$, let us show by induction that, for all $q$ in $\mathbb{N}^*$ such that $2q \leq n$, we have

$$(3.5) \qquad x^q P_n(x) = \beta^0_{n,q} P_{n-2q}(x) + \sum_{n-2q+1 \leq s \leq n+q} \beta^s_{n,q} P_s(x),$$

where $\beta^0_{n,q} = 0$ if $n = 2q$.

For $q = 1$, it is nothing but equality (3.4).

Suppose now that (3.5) is true for $q$ and that $2(q+1) \leq n$. Then we have

$$x^{q+1} P_n(x) = x(x^q P_n(x))$$

$$= \beta^0_{n,q} x P_{n-2q}(x) + \sum_{n-2q+1 \leq s \leq n+q} \beta^s_{n,q} x P_s(x)$$

$$= \beta^0_{n,q} \left\{ \beta^0_{n,n-2q} P_{n-2q-2}(x) + \sum_{n-2q-1 \leq s' \leq n-2q+1} \beta^{s'}_{n-2q,1} P_{s'}(x) \right\}$$

$$+ \sum_{n-2q+1 \leq s \leq n+q} \beta^s_{n,q} \left\{ \beta^0_{s,1} P_{s-2}(x) + \sum_{s-1 \leq s'' \leq s+1} \beta^{s''}_{s,q} P_{s''}(x) \right\}$$

$$= a_3 \beta^0_{n,q} A^3_{n-2q} P_{n-2(q+1)}(x) + \sum_{n-2(q+1)+1 \leq s \leq n+q+1} \beta^s_{n,q+1} P_s(x).$$

Hence there exist $(\beta^s_{n,q+1})$ such that

$$x^{q+1} P_n(x) = \beta^0_{n,q+1} P_{n-2q}(x) + \sum_{n-2q+1 \leq s \leq n+q} \beta^s_{n,q} P_s(x),$$

where $\beta^0_{n,q+1} = a_3 \beta^0_{n,q} A^3_{n-2q}$ and $\beta^0_{n,q+1} = 0$ if $n = 2(q+1)$.

*Proof of* (c). Since $P_1(x) = \frac{1}{V_F(m_o)}(m - m_o)$, it is easy to show by induction that

$$P_n(x) = \frac{1}{(V_F(m_o))^n} x^n + \sum_{0 \leq q \leq n-1} \alpha_q x^q,$$

and this concludes the proof. □

**4. Characterization of the cubic families in the Meixner way.** This section is devoted to the characterization of the 2-orthogonal polynomials on $\mathbb{R}$ with exponential generating function.

We say that the generating function of the sequence of polynomials $Q_n$ is exponential if there exist $r > 0$ and two real analytic functions $a$ and $b$ defined on $]-r, r[$ such that, for all $z$ in $]-r, r[$,

$$(4.1) \qquad \sum_{n \in \mathbb{N}} Q_n(x) \frac{z^n}{n!} = \exp\{a(z)x + b(z)\}.$$

Families $(Q_n/n!)$ satisfying (4.1) are considered in the literature under the name of Sheffer polynomials of type 0; see Sheffer (1939) and Rainville [(1960), Chapter 13]. Actually they are slightly more general since convergence in a neighborhood of 0 is not required and (4.1) is considered as an identity between formal series in $z$.

THEOREM 4.1. *Let $F$ be a NEF on $\mathbb{R}$ and let $\mu$ be an element of $F$ with mean $m_o$. Suppose that $(Q_n)_{n\in\mathbb{N}}$ is a family of $\mu-2$-orthogonal polynomials such that $Q_n$ is of degree $n$. Then the generating function of $(Q_n)_{n\in\mathbb{N}}$ is exponential if and only if there exists $t \in \mathbb{R}^*$ such that, for all $n \in \mathbb{N}$,*

$$Q_n(x) = t^n P_n(x),$$

*where $(P_n)$ is defined by (2.7).*
*In this case, $a(z) = \psi_\mu(tz + m_o)$ and $b(z) = -k_\mu(\psi_\mu(a(z)))$.*

PROOF. Up to $\tilde{Q}_n = Q_n/Q_o$, we can suppose $Q_o = 1$.

$\Leftarrow$ Is obvious.
$\Rightarrow$ There exist $r > 0$ such that, for all $z \in \,]-r, r[$,

$$\int \left( \sum_{n\in\mathbb{N}} Q_n(x) \frac{z^n}{n!} \right) \mu(dx) = \int \left( \sum_{n\in\mathbb{N}} Q_n(x) Q_o(x) \frac{z^n}{n!} \right) \mu(dx)$$
$$= \int Q_o(x)^2 \mu(dx)$$
$$= 1.$$

On the other hand, writing the generating function of $(Q_n)$ as in (4.1), we have

$$\int \left( \sum_{n\in\mathbb{N}} Q_n(x) \frac{z^n}{n!} \right) \mu(dx) = \int \exp\{a(z)x + b(z)\} \mu(dx)$$
$$= \exp\{k_\mu(a(z)) + b(z)\}.$$

Hence

(4.2) $$b(z) = -k_\mu(a(z)).$$

Proceeding similarly, we have that

(4.3) $$\int \left( \sum_{n\in\mathbb{N}} Q_n(x) Q_1(x) \frac{z^n}{n!} \right) \mu(dx) = \left( \int Q_1(x)^2 \mu(dx) \right) z.$$

Then $Q_1$ is a polynomial of degree 1 in $x$. Therefore there exists $u \in \mathbb{R}^*$ and $v \in \mathbb{R}$ such that

(4.4) $$Q_1 = ux + v.$$

Since $\int Q_1(x)Q_o(x)\mu(dx) = \int Q_1(x)\mu(dx) = 0$, we get $b = -um_o$ and

$$\int Q_1(x)^2 \mu(dx) = \int u^2(x - m_o)(x - m_o)\mu(dx) = u^2 V_F(m_o).$$

Furthermore, using (4.2)–(4.4), we get

$$\left(\int Q_1(x)^2 \mu(dx)\right) z = \int \exp\{a(z)x + b(z)\}Q_1(x)\mu(dx)$$

$$= u\int (x - m_o) \exp\{a(z) - k_\mu(a(z))\}\mu(dx)$$

$$= u\int (x - m_o) P(a(z), \mu)(dx)$$

$$= u[k'_\mu(a(z)) - m_o],$$

and we deduce that

$$u^2 V_F(m_o)z = u[k'_\mu(a(z)) - m_o].$$

Therefore, $k'_\mu(a(z)) = u V_F(m_o) z + m_o$, that is, $a(z) = \psi_\mu(u V_F(m_o) z + m_o)$ and $t = u V_F(m_o)$.

Finally, we obtain

$$\sum_{n \in \mathbb{N}} Q_n(x) \frac{z^n}{n!} = f_\mu(x, V_F(m_o)uz + m_o). \qquad \square$$

COROLLARY 4.1. *Let $F$ be a NEF on $\mathbb{R}$ and let $\mu$ be an element of $F$ with mean $m_o$. Then there exists a family of $\mu - 2$-orthogonal polynomials with an exponential generating function if and only if $F$ is cubic.*

PROOF. Follows easily from Theorems 3.1. and 4.1. $\square$

**5. Proof of Theorem 2.1.** This section is devoted to the proof of Theorem 2.1.

Let $x$ be a fixed real number. The function $m \mapsto f_\mu(x, m)$ is real analytic in the interval $M_F$. Therefore, there exist an open set $\Omega$ of $\mathbb{C}$ containing $M_F$, an open set $U$ of $\mathbb{C}$ containing $\Theta(\mu)$ and analytic functions $\psi_1$ and $k_1$ such that $\psi_1|_{M_F} = \psi_\mu$, $k_1|_{\Theta(\mu)} = k_\mu$ and $\psi_1(\Omega) \subset U$. For all $x \in \mathbb{R}$, the function $z \mapsto f_1(x, z) = \exp\{x\psi_1(z) - k_1(\psi_1(z))\}$ is analytic on $\Omega$.

Since $0 \in \Theta(\mu)$ and $m_o \in M_F$, there exist $\alpha > 0$ and $r > 0$ such that:

(i) $[-\alpha, \alpha] \subset \Theta(\mu)$;
(ii) the open disk $D(m_o, r) \subset \Omega$;
(iii) $\psi_1(D(m_o, r)) \subset D(0, \frac{\alpha}{3})$;
(iv) $V_1 = (\psi'_1)^{-1}$ has no zero in $D(m_o, r)$.

Let us show that the function $\phi(z) = f_1(\cdot, z)$ is well defined and continuously differentiable from $D(m_o, r)$ into $L^2(\mu)$.

Since
$$\int_{-\infty}^{+\infty} |f_1(x,z)|^2 \mu(dx)$$
$$= \int_{-\infty}^{+\infty} \exp\{x(\psi_1(z) + \overline{\psi_1(z)}) - (k_1(\psi_1(z)) + \overline{k_1(\psi_1(z))})\}\mu(dx)$$

and $\theta = \psi_1(z) + \overline{\psi_1(z)} \in [\frac{-2\alpha}{3}, \frac{2\alpha}{3}] \subset \Theta(\mu)$, then
$$\int_{-\infty}^{+\infty} |f_1(x,z)|^2 \mu(dx) < +\infty,$$

and $\phi$ is well defined.

For the differentiability, if $z_o$ is in $D(m_o, r)$, then $(\frac{\partial}{\partial z} f_1)(\cdot, z_o)$ is an element of $L^2(\mu)$ because
$$\int_{-\infty}^{+\infty} \left|\frac{x-z}{V_1(z)}\right|^2 \exp\{x(\psi_1(z) + \overline{\psi_1(z)})\}\mu(dx) < +\infty.$$

We will verify that $\frac{f_1(\cdot, z_o + h) - f_1(\cdot, z_o)}{h} - \frac{\partial}{\partial z} f_1(\cdot, z)$ converges to 0 in $L^2(\mu)$ when $h$ converges to 0, that is,

(5.1) $\quad \lim_{h \to 0} \frac{1}{h^2} \int_{-\infty}^{+\infty} \left|f_1(x, z_o + h) - f_1(x, z_o) - h\left(\frac{\partial}{\partial z} f_1\right)(x, z_o)\right|^2 \mu(dx) = 0.$

Writing the Taylor formula with integral remainder,
$$f(h) - f(0) - hf'(0) = h^2 \int_0^1 (1-u) f''(z_o + uh) \, du,$$

(5.1) is equivalent to

(5.2) $\quad \lim_{h \to 0} \frac{1}{h^2} \int_{-\infty}^{+\infty} h^4 \left|\int_0^1 \frac{\partial^2}{\partial z^2} f_1(x, z_o + uh)(1-u)\,du\right|^2 \mu(dx) = 0.$

But we have that
$$\frac{\partial^2}{\partial z^2} f_1(x, z) = \left[\left(\frac{x-z}{V_1(z)}\right)^2 - \frac{V_1(z) + (x-z)V_1'(z)}{V_1^2(z)}\right] f_1(x, z).$$

Then $\frac{\partial^2}{\partial z^2} f_1(x, z_o + h)/f(x, z_o + h)$ is a second-degree polynomial in $x$ whose coefficients are continuous in $h$. It is then bounded in all compact $|h| \leq h_o$. Hence
$$\left|\frac{\partial^2}{\partial z^2} f_1(x, z_o + h)\right| \leq |ax^2 + bx + c|e^{\theta_o x} \qquad \text{with } |\theta_o| < \frac{2\alpha}{3}.$$

This with the dominated convergence theorem implies that

$$\lim_{h \to 0} \int_{-\infty}^{+\infty} \left| \int_0^1 (1-u) \frac{\partial^2}{\partial z^2} f_1(x, z_o + uh) \, du \right|^2 \mu(dx)$$
$$= \int_{-\infty}^{+\infty} \left| \int_0^1 (1-u) \frac{\partial^2}{\partial z^2} f_1(x, z_o) \, du \right|^2 \mu(dx) < +\infty,$$

and (5.2) follows.

Hence $\phi$ is continuously differentiable and so it is analytic on $D(m_o, r)$ and, in particular, if

$$b_n = \left( \int_{-\infty}^{+\infty} P_n(x)^2 \mu(dx) \right)^{1/2} = \|P_n(x)\|_{L^2(\mu)},$$

then $\sum_{n \in \mathbb{N}} b_n \frac{(m-m_o)^n}{n!}$ has nonzero radius of convergence.

Faculté des Sciences
Université de Sfax
B.P. 802
Tunisie
e-mail: abdelhamid.hassairi@fss.rnu.tn